\documentclass[a4paper,10pt]{article}
\usepackage[latin1]{inputenc}
\usepackage[T1]{fontenc}
\usepackage[francais,english]{babel}
\usepackage{amsmath,amscd,hyperref}   
\usepackage{amsfonts}
\usepackage{amssymb}
\usepackage{amsthm}
\usepackage{color}
\usepackage{enumerate}
 \usepackage{verbatim}
 
\newcommand{\md}{\mathrm{d}}
\newcommand{\lag}{\langle}
\newcommand{\rag}{\rangle}

\newcommand{\la}{\lambda}

\newcommand{\bo}{\boldsymbol}

\newcommand{\HH}{\mathcal{H}}

\newcommand{\KK}{\mathcal{K}}
\newcommand{\LL}{\mathcal{L}}
\newcommand{\MM}{\mathcal{M}}

\newcommand{\OO}{\mathcal{O}}
\newcommand{\PP}{\mathcal{P}}
\newcommand{\QQ}{\mathcal{Q}}

\renewcommand{\SS}{\mathcal{S}}

\newcommand{\R}{\mathbb{R}}
\newcommand{\C}{\mathbb{C}}
\newcommand{\N}{\mathbb{N}}

\newtheorem{theorem}{Theorem}[section]
\newtheorem{lemma}[theorem]{Lemma}
\newtheorem{e-proposition}[theorem]{Proposition}

\newtheorem{e-definition}[theorem]{Definition\rm}
\newtheorem{remark}{\it Remark}

\title{Global exact controllability of a \textsc{1d} {S}chr\"odinger equations with a polarizability term}
\author{
Morgan \textsc{Morancey}\footnote{CMLS UMR 7640, Ecole Polytechnique, 91128 Palaiseau, France; CMLA ENS Cachan, 61 av du Pr\'esident Wilson, 94235 Cachan; 
email: Morgan.Morancey@cmla.ens-cachan.fr},
Vahagn \textsc{Nersesyan}\footnote{ Laboratoire de Math\'ematiques, UMR CNRS 8100, Universit\'e de Versailles-Saint-Quentin-Yvelines, F-78035 Versailles, France;  e-mail: Vahagn.Nersesyan@math.uvsq.fr}
\thanks{The authors were partially supported by the ANR grants  EMAQS and STOSYMAP No. ANR-2011-BS01-017-01 and  ANR-2011-BS01015-01}
}
\date{}

\begin{document}

\maketitle

\begin{abstract}
\selectlanguage{english}
We consider a quantum particle in a $1$\textsc{d} interval submitted to a potential. The evolution of this particle is controlled using an external electric field. Taking into account the so-called polarizability term in the model (quadratic with respect to the control),  we prove global exact controllability in a suitable space for arbitrary potential and arbitrary dipole moment. This term is relevant both from the mathematical and physical points of view. The proof uses tools from the bilinear setting and a perturbation argument.

\vskip 0.5\baselineskip

\selectlanguage{francais}
\begin{center} {\bf R\'esum\'e}  \end{center}

\noindent
On consid\`ere une particule quantique dans un intervalle $1$\textsc{d}, soumise \`a un potentiel. L'\'evolution de cette particule est contr\^ol\'ee par un champ \'electrique ext\'erieur. En prenant en compte dans le mod\`ele le terme  dit de polarisabilit\'e (quadratique par rapport au contr\^ole), on prouve la contr\^olabilit\'e exacte globale dans un espace appropri\'e pour des potentiels et des moments dipolaires arbitraires. Ce terme est int\'eressant \`a la fois d'un point de vue math\'ematique et physique. La preuve utilise des outils issus du cadre bilin\'eaire et un argument de perturbation.

\end{abstract}

\selectlanguage{francais}
\section*{Version fran\c{c}aise abr\'eg\'ee}

On consid\`ere une particule quantique unidimensionnelle soumise \`a l'action d'un potentiel $V$. La particule est repr\'esent\'ee par sa fonction d'onde $\psi$ dont l'\'evolution est contr\^ol\'ee par un champ \'electrique ext\'erieur d'amplitude r\'eelle $u$.  En notant $\mu_1$   le moment dipolaire et $\mu_2$ le moment de polarisabilit\'e,  l'\'evolution de la fonction d'onde est donn\'ee par le syst\`eme de Schr\"odinger avec polarisabilit\'e
\begin{equation} \label{syst_pola}
\left\{
\begin{aligned}
&\begin{aligned}
i \partial_t \psi = &\left( -\partial^2_{xx} + V(x) \right) \psi -u(t) \mu_1(x) \psi 
 - u(t)^2 \mu_2(x) \psi,
\end{aligned}
&(t,x)& \in (0,T)\times (0,1),
\\
&\psi(t,0) = \psi(t,1) = 0,  \\
&\psi(0,x) = \psi_0(x).  
\end{aligned}
\right.
\end{equation}
Si la prise en compte du terme de polarisabilit\'e est int\'eressante du point de vue physique (par exemple dans le cas de contr\^oles de fortes amplitudes~\cite{Dion_2}), du point de vue math\'ematique, ce terme a permis de montrer la contr\^olabilit\'e dans des cas o\`u le moment dipolaire est insuffisant pour conclure (voir par exemple~\cite{CGLT}, \cite{Morancey_polarisabilite}, \cite{BCC_polarisabilite}).

Pour $V \in L^2((0,1),\R)$, on note $\lambda_{k,V}$ et $\varphi_{k,V}$ les valeurs propres (en ordre croissant) et vecteurs propres de l'op\'erateur $A_V$ d\'efini sur le domaine $D(A_V) := H^2 \cap H^1_0((0,1),\C)$ par $A_V \psi := \left( -\partial^2_{xx} + V(x) \right) \psi$. 
On d\'efinit les \'etats propres par
$
\Phi_{k,V}(t,x) := e^{-i \lambda_{k,V} t}  \varphi_{k,V}(x),  $ $ (t,x) \in [0,+\infty) \times (0,1), \: k \in \N^*.
$
Pour $s>0$, l'espace $H^s_{(V)}:=D(A_V^{s/2})$ est muni de la norme
$\| \psi \|_{H^s_{(V)}} := \left( \sum_{k=1}^{+\infty} | k^s \lag \psi , \varphi_{k,V} \rag |^2 \right)^{\frac{1}{2}}.$
On note $\SS$ la sph\`ere unit\'e de $L^2((0,1),\C)$.

Dans le cadre bilin\'eaire (c'est-\`a-dire pour le syt\`eme (\ref{syst_pola}) avec $\mu_2=0$), en combinant les r\'esultats de contr\^ole exact local, dans $H^3_{(0)}$, autour de $\Phi_{1,0}$ de Beauchard et Laurent~\cite{BeauchardLaurent} et la contr\^olabilit\'e approch\'ee de $\varphi_{1,0}$ dans $H^3$ du second auteur~\cite{Nersesyan10}, on obtient la contr\^olabilit\'e exacte globale dans $\SS \cap H^{3+\epsilon}_{(0)}$   pour~$V=0$ sous des hypoth\`eses favorables sur $\mu_1$. Ces deux r\'esultats sont principalement bas\'es sur l'\'etude de lin\'earis\'es du syst\`eme   au voisinage de trajectoires associ\'ees au contr\^ole nul. En utilisant le fait que (au moins formellement) le syst\`eme~(\ref{syst_pola}) avec $\mu_2\in L^2((0,1),\R)$ quelconque a    le m\^eme   lin\'earis\'e au voisinage de telles trajectoires  que dans le cas bilin\'eaire $\mu_2=0 $, conjointement \`a un argument de perturbation utilis\'e par les auteurs dans~\cite{MoranceyNersesyan13} dans le cadre du contr\^ole simultann\'e de syst\`emes bilin\'eaires, on prouve le r\'esultat suivant
\begin{theorem} \label{th_global_exact_fr}
Pour tout $V, \mu_1 \in H^6((0,1),\R)$ le syst\`eme (\ref{syst_pola}) est globalement exactement contr\^olable dans~$H^6_{(V)}$, g\'en\'eriquement par rapport \`a $\mu_2 \in H^6((0,1),\R)$. \end{theorem}

Par rapport au mod\`ele bilin\'eaire, la prise en compte du terme de polarisabilit\'e, permet de conclure \`a la contr\^olabilit\'e dans des cas o\`u la contr\^olabilit\'e \'etait fausse ou ouverte (par exemple $V$ arbitraire et $\mu_1=0$ ou $\mu_1 \notin \QQ_V$ comme d\'efini dans~\cite{MoranceyNersesyan13}).

\selectlanguage{english}
\section{Introduction}
\label{}

We consider the evolution of a $1$\textsc{d} quantum particle given by (\ref{syst_pola}).
The real valued functions $V,\mu_1,$ and $\mu_2$, respectively, the potential, the dipole moment, and the polarizability moment, are given. The control~$u(t)$ is real valued. The following theorem is the main result of this paper.   
\begin{theorem} \label{th_global_exact}
For any  $V, \mu_1 \in H^6((0,1),\R)$, system (\ref{syst_pola}) is globally exactly controllable in~$H^6_{(V)}$ generically with respect to  $\mu_2 \in H^6((0,1),\R)$. More precisely,  there is a residual set   $\QQ_{V,\mu_1}$ in $H^6((0,1),\R)$ such that if  $\mu_2 \in \QQ_{V,\mu_1}$, then for any $\psi_0, \psi_f \in \SS \cap H^6_{(V)}$, there is  $T>0$ and $u \in H^1_0((0,T),\R)$ such that the solution of (\ref{syst_pola}) satisfies $\psi(T) = \psi_f$.
\end{theorem}

Essentially with the same proof one can establish the same exact controllability property in the case where the term $u(t)^2 \mu_2(x) \psi$ in (\ref{syst_pola}) is replaced by a higher degree term $\sum_{j=2}^m u^j\mu_j \psi.$ We choose $m=2$ for the sake of simplicity of presentation.

\medskip 
\textit{Review of previous results.} The controllability properties of quantum particles were first studied for the bilinear model (i.e., for~(\ref{syst_pola}) with $\mu_2=0$). In~\cite{BeauchardLaurent}, Beauchard and Laurent proved local exact controllability in $H^3_{(0)}$ around $\Phi_{1,0}$ by studying the controllability of the linearized system around the trajectory $(u \equiv 0, \Phi_{1,0})$. The simultaneous global exact controllability of an arbitrary (finite) number of such bilinear equations was studied by the authors in~\cite{MoranceyNersesyan13} for arbitrary potentials.
For multidimensional domains, we mention the simultaneous approximate controllability property obtained by Boscain, Caponigro, Chambrion, Mason, Sigalotti \cite{CMSB09}, \cite{CBCS11} through geometric techniques
based on the   exact controllability of the Galerkin approximations.
The approximate controllability in Sobolev spaces towards the state  $\varphi_{1,V}$  is obtained by the second author using Lyapunov techniques~\cite{Nersesyan10}.     The first  controllability results for systems having a   polarizability term are established for  finite-dimensional models and are due to  Coron, Grigoriu, Lefter, and Turinici~\cite{Turinici07}, \cite{GrigoriuLefterTurinici09}, \cite{CGLT}, \cite{Grigoriu12}. They proved exact controllability under the same assumptions as for the bilinear model. They also proved stabilization of the first eigenstate using either discontinuous feedback laws or time oscillating periodic feedback laws in a setting where the dipole moment was not sufficient to conclude. The strategy based on time oscillating feedback laws has been extended to the infinite dimensional model~(\ref{syst_pola}) by the first author~\cite{Morancey_polarisabilite}. Finally, geometric technics  were   applied to the polarizability system by Boussaid, Caponigro, Chambrion in~\cite{BCC_polarisabilite} leading to global approximate controllability.

\medskip
\textit{Structure of the article.} 
First, we prove in Section~\ref{sect_approche}   approximate controllability towards the ground state, adapting Lyapunov arguments from the bilinear setting. Still using tools from the bilinear setting, we prove in  Section~\ref{sect_exact}  local exact controllability around the ground state, in $H^5_{(V)}$, with controls in~$H^1_0((0,T),\R)$.  Gathering these results, we get  global exact controllability under favourable hypotheses on~$V$ and $\mu_1$ in Section~\ref{sect_global_exact}. Then, we conclude the proof using a perturbation argument.

\section{Approximate controllability towards the ground state}
\label{sect_approche}

  From the arguments of the proofs of   \cite[Proposition 3.1]{NersesyanNersisyan1D}  and \cite[Proposition 5]{BeauchardLaurent} it easily follows    that for any $V, \mu_1, \mu_2 \in H^5((0,1),\R)$, $T>0$, $\psi_0 \in H^5_{(V)}$, and $u \in H^1_0((0,T),\R)$, system  (\ref{syst_pola}) has a unique weak solution $\psi \in C([0,T],H^5)\cap  C^1([0,T],H^3_{(V)})$. Furthermore,   the mapping which sends $(\psi_0, u)$ to the solution $\psi$ is~$C^1$.
As  $u(T)=0$, it comes that $\psi(T) \in H^5_{(V)}$. When  $\psi_0 \in H^6_{(V-u(0)\mu_1-u(0)^2 \mu_2)}$,   $u \in C^2([0,T],\R)$, and $\dot{u}(0)=0$,  the   solution $\psi $ belongs to $C([0,T],H^6)$. Moreover, if $\dot{u}(T)=0$, then $\psi(T,\psi_0,u) \in H^6_{(V-u(T)\mu_1 -u(T)^2 \mu_2)}$.

\medskip
Let us introduce the following Lyapunov function
\begin{equation} \label{def_lyapu}
\LL(z) := \gamma \| \left( - \partial^2_{xx} + V \right)^3 \PP z \|_{L^2}^2 + 1 - |\lag z, \varphi_{1,V}\rag|^2, \quad z \in \SS \cap H^6_{(V)},
\end{equation}
where $\PP$ is the orthogonal projection in $L^2$ onto the closure of the vector space spanned by $\{ \varphi_{k,V} \, ; \, k \geq 2 \}$ and $\gamma>0$ is a constant which will be precised later on. 
This Lyapunov function has already been used in the bilinear setting by the second author and Beauchard in~\cite{Nersesyan10}, \cite{BeauchardNersesyan} and adapted to study simultaneous controllability in~\cite{MoranceyNersesyan13} by the authors.
We assume that the functions $V,\mu_1 \in H^6((0,1),\R)$ are such that 
\smallskip
\begin{itemize}
\item[$\boldsymbol{\mathrm{(C_1)}}$] $\lag \mu_1 \varphi_{1,V} , \varphi_{k,V} \rag \neq 0$, for all $k \in \N^*$,
\smallskip
\item[$\boldsymbol{\mathrm{(C_2)}}$] $\lambda_{1,V} - \lambda_{j,V} \neq \lambda_{p,V} - \lambda_{q,V}$, for all $j,p,q \geq 1$ and $j \neq 1$.
\smallskip
\end{itemize}

\begin{theorem} \label{th_global_approche}
Let $V, \mu_1, \mu_2 \in H^6((0,1),\R)$ be such that Conditions $\bo{(\mathrm{C_1})}$ and $\bo{(\mathrm{C_2})}$ are satisfied. For any $\psi_0 \in \SS \cap H^6_{(V)}$ satisfying $\lag \psi_0 ,\varphi_{1,V} \rag \neq 0$ and for any $\varepsilon>0$, there are $T>0$ and $u \in C^2_0((0,T),\R)$ such that
$\| \psi(T, \psi_0 ,u) - \varphi_{1,V} \|_{H^5} < \varepsilon.
$\end{theorem}

\noindent
\textit{Proof.} For $\mu_2 = 0$, the proof is given in~\cite[Theorem 2.3]{Nersesyan10}. The adaptation to $\mu_2 \in H^6((0,1),\R)$ is straightforward, we only recall the scheme of the proof. 
Since for every $\psi_0 \in H^6_{(V)}$, the linearization   of~(\ref{syst_pola}) around the trajectory $\psi(\cdot, \psi_0, 0)$ is the same for $\mu_2=0$ or any other $\mu_2\in H^6((0,1),\R)$, we deduce from~\cite[Proposition 2.6]{Nersesyan10} the following~lemma.  
\begin{lemma} \label{lemme_decroissance_lyapu}
Let $V, \mu_1, \mu_2 \in H^6((0,1),\R)$ be such that Conditions $\bo{(\mathrm{C_1})}$ and  $\bo{(\mathrm{C_2})}$ are satisfied.  For  any   $\psi_0 \in \SS \cap H^6_{(V)}$ satisfying $\lag \psi_0, \varphi_{1,V} \rag \neq 0$ and $\LL(\psi_0) > 0$, there exist a time $T>0$ and a control $u \in C^2_0((0,T),\R)$ such that
$\LL(\psi(T,\psi_0,u)) < \LL(\psi_0).
$\end{lemma}

Let us take any  $\psi_0 \in \SS \cap H^6_{(V)}$ satisfying $\lag \psi_0, \varphi_{1,V} \rag \neq 0$  and let us choose the constant $\gamma >0$ in (\ref{def_lyapu})    such that $\LL(\psi_0) < 1$. If $\LL(\psi_0) > 0$, we define
\begin{gather*}
\KK := \Big\{ \psi_f \in H^6_{(V)} \, ; \, \psi(T_n,\psi_0,u_n) \underset{n \to \infty}{\longrightarrow} \psi_f, \text{ in } H^5 \text{ where } T_n >0,
 u_n \in C^2_0((0,T_n),\R) \Big\}.
\end{gather*}
The infimum of $\LL$ on $\KK$ is attained, i.e., there is $e \in \KK$ such that 
$\LL(e) = \inf_{\psi \in \KK} \LL(\psi).
$ This gives that $\LL(e) \leq \LL(\psi_0) < 1$, hence $\lag e , \varphi_{1,V} \rag \neq 0$. Using   Lemma~\ref{lemme_decroissance_lyapu}, it comes that if $\LL(e) >0$ then there are $T>0$ and $u \in C^2_0((0,T),\R)$ such that $\LL(\psi(T,e,u)) < \LL(e)$. As $\psi(T,e,u) \in \KK$, this contradicts the definition of~$e$. Then, $\LL(e) = 0$. This leads to $\varphi_{1,V} \in \KK$ and concludes the proof of Theorem~\ref{th_global_approche}.

\hfill $\blacksquare$

\section{Local exact controllability around the ground state}
\label{sect_exact}

In this section, we     prove local exact controllability  around the ground state $\varphi_{1,V}$  in $H^5_{(V)}$  with controls in $H^1_0((0,T),\R)$. We assume that the functions $V,\mu_1 \in H^5((0,1),\R)$ satisfy
\smallskip
\begin{itemize}
\item[$\boldsymbol{\mathrm{(C_3)}}$] there exists $C>0$ such that
$\left| \lag \mu_1 \varphi_{1,V} , \varphi_{k,V} \rag \right| \geq \frac{C}{k^3}, \text{ for all } k \in \N^*.
$\smallskip
\end{itemize}

\begin{theorem} \label{th_local_exact}
Let $V, \mu_1, \mu_2 \in H^5((0,1),\R)$ be such that Condition $\bo{(\mathrm{C_3})}$ is satisfied. Let $T>0$. There exist $\delta>0$ and a $C^1$ map
$\Gamma : \OO_T \longrightarrow H^1_0((0,T),\R)
$, where
$\OO_T := \left\{ \psi_f \in \SS \cap H^5_{(V)} \, ;\, \| \psi_f - \Phi_{1,V}(T)\|_{H^5} < \delta \right\},
$ such that $\Gamma(\Phi_{1,V}(T))=0$, and for any $\psi_f \in \OO_T$, the solution of (\ref{syst_pola}) with initial condition $\psi_0=\varphi_{1,V}$ and control $u=\Gamma(\psi_f)$ satisfies $\psi(T) = \psi_f$.
\end{theorem}

\noindent
\textit{Proof.}   In the case where $\mu_2 = 0$ and $V=0$, the proof is exactly the one of \cite[Theorem 2]{BeauchardLaurent}. Let  
$\HH := \left\{ \psi \in H^5_{(V)} \, ; \, \Re(\lag \psi , \varphi_{1,V} \rag ) = 0 \right\}
$, and let $\PP_{\HH}$ be the orthogonal projection in $L^2((0,1),\C)$ onto $\HH$.  Then the  end-point map $\Theta_{T}:   u \in  H^1_0((0,T),\R)   \mapsto  \PP_{\HH} \big( \psi(T, \varphi_{1,V},u) \big) \in  \HH$,  is $C^1$ and  its differential at $0$ is given by $\md \Theta_T(0).v = \Psi(T)$, where $\Psi$ is the solution of  
\begin{equation}\label{Lin}
 i \partial_t \Psi = \left( -\partial^2_{xx} + V(x) \right) \Psi - v(t) \mu_1(x) \Phi_{1,V},  \quad (t,x)  \in (0,T) \times (0,1) 
\end{equation} 
 with homogeneous Dirichlet boundary conditions and   $\Psi(0,x) = 0$. 	 
Rewriting this   in the Duhamel form, we get 
$$
\Psi(T) = i \sum_{k=1}^{+\infty} \lag \mu_1 \varphi_{1,V}, \varphi_{k,V} \rag \int_0^T v(t) e^{ i( \lambda_{k,V} - \lambda_{1,V}) t } \md t \, \Phi_{k,V}(T).
$$
Using Condition $\bo{(\mathrm{C_3})}$, the asymptotics of eigenvalues $\la_{k,V}$,  and \cite[Corollary 2]{BeauchardLaurent}, we get the existence of a continuous linear map
$\MM :  \HH \mapsto  L^2((0,T),\R)
$ such that for any $\Psi_f \in \HH$, the function $w := \MM(\Psi_f)$ solves the following moment problem
\begin{equation*}  
\left\{
\begin{aligned}
&\int_0^T w(t) \md t = 0, \quad 
 \int_0^T w(t) (T-t) \md t = \frac{1}{\lag \mu_1 \varphi_{1,V}, \varphi_{1,V} \rag} \lag \Psi_f , \Phi_{1,V}(T) \rag,
\\
&\int_0^T w(t) e^{i (\lambda_{k,V}- \lambda_{1,V})t} \md t = \frac{\lambda_{1,V} - \lambda_{k,V}}{\lag \mu_1 \varphi_{1,V}, \varphi_{k,V} \rag}  \lag \Psi_f, \Phi_{k,V}(T) \rag, \quad \forall k \geq 2.
\end{aligned}
\right.
\end{equation*} 
Then the mapping $\Psi_f \in \HH \mapsto \int_0^t w(\tau) \md \tau \in H^1_0((0,T),\R)$,
is a continuous right inverse for the differential~$\md \Theta_T(0)$.
Finally, applying the inverse mapping theorem to $\Theta_T$ at $u=0$, using the conservation of the $L^2$ norm and the hypothesis that $\psi_f \in \SS$, we complete the    proof of Theorem~\ref{th_local_exact}.

\hfill $\blacksquare$

\begin{remark} Let us notice that the linearized system (\ref{Lin}) is controllable in $H^3_{(V)}$ with controls in $L^2((0,T),\R)  $, but we cannot conclude to controllability for (\ref{syst_pola}), since we do not know if the latter is well posed in $H^3_{(V)}$ for $u\in L^2((0,T),\R)  $. Indeed, in that case  $u^2$ will be in $L^1((0,T),\R),$ which does not allow to apply the results of \cite{BeauchardLaurent}. 
\end{remark}

\section{Global exact controllability}
\label{sect_global_exact}

Combining the properties of global approximate controllability    obtained in Theorem~\ref{th_global_approche} and   local exact controllability obtained in Theorem~\ref{th_local_exact}, we obtain global exact controllability of (\ref{syst_pola}) in $H^6_{(V)}$ under favourable hypotheses on $V$ and $\mu_1$.
\begin{theorem} \label{th_global_exact_partiel}
Let $V, \mu_1, \mu_2 \in H^6((0,1),\R)$ be such that Conditions $\bo{(\mathrm{C_2})}$ and $\bo{(\mathrm{C_3})}$ are satisfied. For any $\psi_0, \psi_f \in \SS \cap H^6_{(V)}$, there is a time $T>0$ and a control $u \in H^1_0((0,T),\R)$ such that the associated solution of (\ref{syst_pola}) satisfies  $\psi(T) = \psi_f$.
\end{theorem}

\noindent
\textit{Proof.}  \textit{First step.} Let $\psi_0, \psi_f \in \SS \cap H^6_{(V)}$ be such that $\lag \psi_0 , \varphi_{1,V}\rag \neq 0$ and $\lag \psi_f ,\varphi_{1,V} \rag \neq 0$. Let   $T_*>0$ be such that $\Phi_{1,V}(T_*)=\varphi_{1,V}$ and let  $\delta>0$ be the radius of local exact controllability in $H^5_{(V)}$ in time $T_*$ given by Theorem~\ref{th_local_exact}.
  Theorem~\ref{th_global_approche} implies the existence of  times $T_0, T_f >0$ and controls $u_0 \in C^2_0((0,T_0),\R)$,  $u_f \in C^2_0((0,T_f),\R)$ such that
$$
\| \psi(T_0,\psi_0,u_0) - \varphi_{1,V} \|_{H^5} + \| \psi(T_f, \overline{\psi_f},u_f) - \varphi_{1,V} \|_{H^5} < \delta.
$$ By Theorem~\ref{th_local_exact}, there exists $u_* \in H^1_0((0,T_*),\R)$ such that 
$\psi(T_*, \varphi_{1,V}, u_*) =   \overline{\psi(T_0,\psi_0,u_0)}.
$ Time reversibility property of (\ref{syst_pola}) implies that, if we define    $u(t) = u_0(t)$ for $t\in[0,T_0]$ and $u(t+T_0) = u_* (T_*-t)$ for $t \in [0,T_*]$,  then $u \in H^1_0((0,T_0+T_*),\R)$ and
$\psi(T_0+T_*,\psi_0,u)= \varphi_{1,V}.
$  
The same arguments lead to the existence of $\tilde{u} \in H^1_0((0,T_f+T_*),\R)$ such that
$\psi(T_f+T_*, \overline{\psi_f}, \tilde{u}) = \varphi_{1,V}.
$ Taking  $T := T_0 + T_f + 2T_*$ and again applying  the time reversibility, we     find $u \in H^1_0((0,T),\R)$ satisfying
$\psi(T,\psi_0,u) = \psi_f.
$

\medskip
\textit{Second step.} It remains to remove the hypotheses $\lag \psi_0 , \varphi_{1,V}\rag \neq 0$ and $\lag \psi_f ,\varphi_{1,V} \rag \neq 0$.
Using time reversibility, it is sufficient to prove that for any $\psi_0 \in \SS $, there are $T>0$ and $u \in C^2_0((0,T)\R)$ such that   $\lag \psi(T,\psi_0,u), \varphi_{1,V} \rag \neq 0$. 
Let $\widehat{\psi}_0 \in \SS \cap H^6_{(V)}$ be such that $\lag \widehat{\psi}_0 , \varphi_{1,V} \rag \neq 0$ and $\| \psi_0 - \widehat{\psi}_0 \|_{L^2} < \sqrt{2}$. 
From the first step, we get the existence of $T>0$ and $\hat{u} \in H^1_0((0,T),\R)$ such that $\psi(T,\widehat{\psi}_0, \hat{u}) = \varphi_{1,V}$. Then, the conservation of the $L^2$ norm implies 
\begin{equation*}
\| \psi(T,\psi_0,\hat{u}) - \varphi_{1,V} \|_{L^2} = \| \psi_0 - \widehat{\psi}_0 \|_{L^2} < \sqrt{2}.
\end{equation*}
Choosing $u \in C^2_0((0,T),\R)$ sufficiently close to $\hat{u}$ in $L^2((0,T),\R)$, we complete  the proof of  Theorem~\ref{th_global_exact_partiel}.

\hfill $\blacksquare$

\medskip
Finally, adapting a perturbation argument from \cite{MoranceyNersesyan13} and applying  Theorem~\ref{th_global_exact_partiel}, we prove   Theorem~\ref{th_global_exact}.

\noindent
 \textit{Proof of Theorem~\ref{th_global_exact}.}   Let $V, \mu_1 \in H^6((0,1),\R)$ and let $\QQ_{V,\mu_1}$ be the set of functions $\mu_2 \in H^6((0,1),\R)$ such that Conditions $\bo{(\mathrm{C_2})}$ and $\bo{(\mathrm{C_3})}$ are satisfied with the functions $V$ and $\mu_1$ replaced, respectively, by $V-2\mu_1 - 4\mu_2$ and $\mu_1+4\mu_2$, i.e., $\QQ_{V,\mu_1} := \left\{ \mu_2 \in H^6((0,1),\R) \, ; \, \text{Conditions } \bo{(\mathrm{C_2'})} \text{ and } \bo{(\mathrm{C_3'})  } \text{ hold}\right\}
$,  where 
\begin{itemize}
\smallskip
\item[$\boldsymbol{\mathrm{(C_2')}}$] $\lambda_{1,V-2\mu_1 - 4\mu_2} - \lambda_{j,V-2\mu_1 - 4\mu_2} \neq \lambda_{p,V-2\mu_1 - 4\mu_2} - \lambda_{q,V-2\mu_1 - 4\mu_2}$, for $j,p,q \geq 1$ and $j \neq 1$.
\smallskip
\item[$\boldsymbol{\mathrm{(C_3')}}$] there exists $C>0$ such that
$ \left| \lag (\mu_1+4\mu_2) \varphi_{1,V-2\mu_1 - 4\mu_2} , \varphi_{k,V-2\mu_1 - 4\mu_2} \rag \right| \geq \frac{C}{k^3}, \text{ for every } k \in \N^*.
$ \end{itemize}

\medskip 
\textit{First step : global exact controllability when $\mu_2 \in \QQ_{V,\mu_1}$.} Let us consider the   equation 
\begin{equation} \label{syst_pola_generique}
i \partial_t \psi =  \left( -\partial^2_{xx} + V(x) - 2 \mu_1(x) - 4 \mu_2(x) \right) \psi 
- \tilde{u}(t) (\mu_1 + 4 \mu_2)(x) \psi  
 - \tilde{u}(t)^2 \mu_2(x) \psi
\end{equation}
  for $(t,x) \in (0,T) \times (0,1)$ and homogeneous Dirichlet boundary conditions.
We denote by $\widetilde{\psi}(T,\psi_0,u)$ its solution at time $T$. Then
\begin{equation} \label{lien_pola_pola_generique}
\widetilde{\psi}(t,\psi_0,u) = \psi(t,\psi_0,u+2) \quad \text{for } t \in [0,T].
\end{equation}
Let $\psi_0, \psi_f \in \SS \cap H^6_{(V)}$ and let $u_1 \in C^2([0,1],\R)$ be such that $u_1(0)=\dot{u}_1(0)=\dot{u}_1(1)=0$ and $u_1(1)=2$. Then 
$
\widetilde{\psi}_0  := \psi(1,\psi_0,u_1) \in \SS \cap H^6_{(V-2\mu_1-4\mu_2)}$ 
 and  $ \overline{\widetilde{\psi}_f}  := \psi(1,\overline{\psi_f},u_1) \in \SS \cap H^6_{(V-2\mu_1-4\mu_2)}$.
As $\mu_2 \in \QQ_{V,\mu_1}$, Theorem~\ref{th_global_exact_partiel} implies the existence of $\tilde{T}>0$ and $\tilde{u} \in H^1_0((0,\tilde{T}),\R)$ such that
$\widetilde{\psi}(\tilde{T},\widetilde{\psi}_0,\tilde{u})= \widetilde{\psi}_f.
$ 
Let $T=2+\tilde{T}$ and $u(t)=u_1(t)$    for   $t  \in [0,1]$, $u(t)=\tilde{u}(t-1)+2$    for   $t  \in  [1,\tilde{T}+1]$, and $u(t)=u_1(1-(t-1-\tilde{T}))  $    for   $t  \in  [\tilde{T}+1, T]$.  
Then time reversibility of (\ref{syst_pola}) and (\ref{lien_pola_pola_generique}) implies
$\psi(T,\psi_0,u) = \psi_f
$ with $u \in H^1_0((0,T),\R)$.

\medskip 
\textit{Second step : genericity.} We conclude the proof of Theorem~\ref{th_global_exact} by showing that $\QQ_{V,\mu_1}$ is residual in $H^6((0,1),\R)$. For any $W \in H^6((0,1),\R)$, let    $\QQ_W$ be the set of    functions $\mu \in H^6((0,1),\R)$ such   that there is $C>0$ satisfying
$\left| \lag \mu \varphi_{1,W+\mu} , \varphi_{k,W+\mu} \rag \right| \geq \frac{C}{k^3}, \quad \forall k \in \N^*
$. 
By \cite[Lemma 5.3]{MoranceyNersesyan13} with $s=6$,   the set $\QQ_W$  is residual in $H^6((0,1),\R)$. Let $W := V-\mu_1 \in H^6((0,1),\R)$. For any $\mu \in \QQ_W$,  if we set $\mu_2 := -\frac{1}{4}( \mu_1 + \mu)$, then $\mu_2\in \QQ_{V,\mu_1}$. This ends the proof of Theorem~\ref{th_global_exact}.

\hfill $\blacksquare$

 \paragraph*{Acknowledgements :}
  The  authors thank K.~Beauchard for fruitful discussions.

\end{document}